\documentclass{amsart}
\usepackage{latexsym}
\usepackage{amsfonts,amssymb,amsmath}
\usepackage{amsthm}
\usepackage{graphicx,color}
\usepackage{amsmath} 
\usepackage{multirow,bigdelim}
\usepackage{cleveref}

\newtheorem{theorem}{Theorem}	
\crefname{theorem}{Theorem}{Theorems}
\newtheorem{lemma}{Lemma}
\newtheorem{corollary}{Corollary}		
	
\crefname{proposition}{Proposition}{Propositions}

\newtheorem{definition}{Definition}

\newtheorem{remark}{Remark}

\crefname{section}{Section}{Sections}
\crefname{theorem}{Theorem}{Theorems}
\crefname{lemma}{Lemma}{Lemmas}
\crefname{corollary}{Corollary}	{Corollaries}			
\crefname{proposition}{Proposition}{Propositions}	
\crefname{claim}{Claim}{Claims}
\crefname{conjecture}{Conjecture}{Conjectures}			
\crefname{definition}{Definition}{Definitions}
\crefname{problem}{Problem}{Problems}
\crefname{example}{Example}{Examples}
\crefname{remark}{Remark}{Remarks}
\crefname{figure}{Figure}{Figures}

\newcommand{\wt}{\widetilde}
\newcommand{\R}{\mathbb{R}}

\DeclareMathOperator{\codim}{codim}

\newfont{\bg}{cmr9 scaled\magstep2}
\newcommand{\bigzerol}{\smash{\lower1.0ex\hbox{\bg 0}}}

\def\<{\color{black}} 
\def\>{\color{black}} 
\makeatother
\begin{document}


\title[Transversality theorems on generic linearly perturbed mappings]{Transversality theorems on generic linearly perturbed mappings}


\author{Shunsuke Ichiki}

\dedicatory{{\it In memory of John Mather }}
\thanks{Research Fellow PD of Japan Society for the Promotion of Science}
\address{
Graduate School of Environment and Information Sciences,  
Yokohama National University, 
Yokohama 240-8501, Japan}
\email{ichiki-shunsuke-jb@ynu.jp}


\subjclass[2010]{57R45, 57R42}

\keywords{generic linear perturbation, transversality, immersion, injection}


\begin{abstract}
In his celebrated paper ``Generic projections'', 
John~Mather has given a striking transversality theorem and 
its applications on generic projections. 
On the other hand, in this paper, two transversality theorems on generic 
linearly perturbed $C^r$ mappings are shown $(r\geq 1)$. 
Moreover, some applications of the two theorems are also given.

\end{abstract}

\maketitle

\section{Introduction}\label{section 1}
Throughout this paper, let $\ell$, $m$ and $n$ stand for positive integers. 
In this paper, unless otherwise stated, 
all manifolds and mappings are assumed to be of class $C^r$ $(r\geq 1)$ 
and all manifolds are assumed to be without boundary and 
to have countable bases. 
%
%

Let $F:U\to \mathbb{R}^\ell$ be a $C^r$ mapping from 
an open subset $U$ of $\R^m$. Then, for any linear mapping 
$\pi:\mathbb{R}^m\to \mathbb{R}^\ell$, set 
%
%
\begin{eqnarray*} 
F_\pi=F+\pi. 
\end{eqnarray*} 
Here, the mapping $\pi$ in $F_\pi=F+\pi$ is restricted to $U$. 

Let $\mathcal{L}(\mathbb{R}^{m},\mathbb{R}^{\ell})$ 
be the space consisting of all linear mappings 
of $\mathbb{R}^{m}$ into $\mathbb{R}^{\ell}$. 
Notice that we have the natural identification 
$\mathcal{L}(\mathbb{R}^{m},\mathbb{R}^{\ell})=(\mathbb{R}^m)^\ell$.  
By $N$, we denote a $C^r$ manifold of dimension $n$. 
For given $C^r$ mappings $f:N \to U$ and $F:U\to \R^\ell$, 
a property of mappings $F_\pi \circ f:N\to \mathbb{R}^\ell$ 
(resp., $\pi \circ f:N\to \mathbb{R}^\ell$) will be said to be 
true for a {\it generic linearly perturbed mapping} 
(resp., a {\it generic projection}) 
if there exists a subset $\Sigma$ with Lebesgue measure zero of 
$\mathcal{L}(\mathbb{R}^{m},\mathbb{R}^{\ell})$ 
such that for any $\pi \in \mathcal{L}(\mathbb{R}^{m},\mathbb{R}^{\ell})-\Sigma$, 
the mapping $F_\pi \circ f:N\to \mathbb{R}^{\ell}$ 
(resp., $\pi \circ f : N\to \mathbb{R}^{\ell}$) has the property. 

In his celebrated paper \cite{GP}, for a given $C^\infty$ embedding 
$f:N\to \R^m$, John~Mather has given 
a striking transversality theorem on 
a generic projection $\pi\circ f:N\to \mathbb{R}^\ell$ $(m>\ell)$, 
where $N$ is a $C^\infty$ manifold 
(for details on this result, see \cite[Theorem~1~(p.~229)]{GP}). 
Moreover, in \cite{GP}, as an application of this result, 
he has also shown that if $f:N\to \R^m$ is a $C^\infty$ embedding and 
$(n,\ell)$ is in the nice range of dimensions 
(for the definition of nice rage of dimensions, refer to \cite{mather6}), 
then a generic projection $\pi\circ f:N\to \mathbb{R}^\ell$ $(m>\ell)$ is stable, 
where $N$ is a compact $C^\infty$ manifold. 

In \cite{GL}, 
an improvement of the transversality theorem of \cite{GP} is given 
by replacing generic projections by generic linear perturbations. 
Namely, in \cite{GL}, for a given $C^\infty$ embedding 
$f:N\to U$ and a given $C^\infty$ mapping 
$F:U\to \R^\ell$, 
a transversality theorem on a generic linearly perturbed mapping 
$F_\pi\circ f:N\to \mathbb{R}^\ell$ is given, 
where $N$ is a $C^\infty$ manifold and 
$\ell$ is an arbitrary positive integer which 
may possibly satisfy $m\leq \ell$. 

Moreover, in \cite{CG}, 
for a given $C^\infty$ immersion or 
a given $C^\infty$ injection 
$f:N\to U$, 
transversality theorems on a generic linearly perturbed mapping 
$F_\pi\circ f:N\to \mathbb{R}^\ell$ are given, 
where $N$ is a $C^\infty$ manifold, $F:U\to \R^\ell$ is a $C^\infty$ mapping and 
$\ell$ is an arbitrary positive integer which 
may possibly satisfy $m\leq \ell$. 

On the other hand, 
in this paper, 
as improvements of some results in \cite{CG}, 
two main transversality theorems (\cref{main,main2} in \cref{section 2}) 
and their applications 
on generic linearly perturbed mapping are given 
in the case where manifolds and mappings are not necessarily of class $C^\infty$. 
%
%
%

The first main theorem 
(\cref{main}) is as follows. 
Let $f:N\to U$ (resp., $F:U\to \R^\ell$) be a $C^r$ immersion 
(resp., a $C^r$ mapping), 
where $N$ is a $C^r$ manifold (for the value of $r$, see \cref{main}). 
Then, 
generally, the composition $F\circ f$ does not necessarily 
yield a mapping transverse 
to the subfiber-bundle of the jet bundle $J^1(N,\mathbb{R}^\ell)$ with 
a fiber $\Sigma^k$, 
where $k$ is a positive integer satisfying $1\leq k \leq \min\{n,\ell\}$ and  
\begin{eqnarray*} 
\Sigma ^k=\left\{j^1g(0)\in J^1(n,\ell)\mid {\rm corank\ }Jg(0)=k\right\}. 
\end{eqnarray*}
Nevertheless, \cref{main} asserts that a 
generic linearly perturbed mapping $F_\pi \circ f$ yields a mapping transverse to 
the subfiber-bundle of $J^1(N,\mathbb{R}^\ell)$ with $\Sigma^k$. 
The second main theorem (\cref{main2}) is a specialized transversality theorem on crossings of 
a generic linearly perturbed mapping $F_\pi \circ f$, 
where $N$ is a $C^r$ manifold, $f:N\to U$ is a given $C^r$ injection and  
$F:U\to \mathbb{R}^\ell$ is a given $C^r$ mapping 
(for the value of $r$, see \cref{main2}). 

For a given $C^2$ immersion (resp., $C^1$ injection) $f:N\to U$ and 
a given $C^2$ mapping (resp., $C^1$ mapping) $F:U\to \R^\ell$, 
the following (1) and (2) (resp., (3)) are obtained as applications of 
\cref{main} (resp., \cref{main2}), where $N$ is a $C^2$ manifold 
(resp., a $C^1$ manifold). 

\begin{enumerate}
\item[(1)] If $(n,\ell)=(n,1)$, 
then a generic linearly perturbed function $F_{\pi}\circ f:N\to \mathbb{R}$ is a Morse function. 
\item[(2)]If $\ell \geq 2n$, 
then a generic linearly perturbed mapping $F_{\pi}\circ f:N\to \mathbb{R}^{\ell}$ is an immersion.  
\item[(3)]If $\ell > 2n$, 
then a generic linearly perturbed mapping $F_{\pi}\circ f:N\to \mathbb{R}^{\ell}$ is an injection. 
\end{enumerate}
Furthermore, by combining the assertions (2) and (3), 
for a given $C^2$ embedding $f:N\to U$ and 
a given $C^2$ mapping $F:U\to \R^\ell$, we get 
the following assertion (4), 
where $N$ is a $C^2$ manifold. 
\begin{enumerate}
\item[(4)]If $\ell > 2n$ and $N$ is compact, 
then a generic linearly perturbed mapping $F_{\pi}\circ f:N\to \mathbb{R}^{\ell}$ is an embedding. 
\end{enumerate}
\par 
\bigskip 
In \cref{section 2}, some definitions are prepared, 
and the two main transversality theorems (\cref{main,main2}) are stated. 
\cref{section 3} (resp., \cref{section 4}) is devoted to
the proof of \cref{main} (resp., \cref{main2}). 
In \cref{section 5}, 
the above assertions (1)--(4) are shown. 
In \cref{app}, the important lemma for the proofs of 
\cref{main,main2} (\cref{abra} in \cref{section 2}) is shown as an appendix. 

\section{Preliminaries and the statements of \cref{main,main2}}\label{section 2}
Firstly, the definition of transversality is given. 
\begin{definition}\label{transverse}
{\rm
Let $N$ and $P$ be $C^r$ manifolds, and $Z$ be a $C^r$ submanifold of $P$ 
($r\geq1$). 
Let $g : N\to P$ be a $C^1$ mapping. 
\begin{enumerate}
\item 
We say that $g:N\to P$ is {\it transverse} to $Z$ 
{\it at $q$} if $g(q)\not\in Z$ or  
in the case of $g(q)\in Z$, the following holds: 
\begin{eqnarray*}
dg_q(T_qN)+T_{g(q)}Z=T_{g(q)}P.
\end{eqnarray*}
\item 
We say that $g:N\to P$ is {\it transverse} to $Z$ 
if for any $q\in N$, 
the mapping $g$ is transverse to $Z$ at $q$. 
\end{enumerate}
}
\end{definition}


For the statement and the proof of \cref{main}, 
some definitions are prepared. 
Let $N$ be a $C^r$ manifold $(r\geq 2)$ and 
$J^1(N,\R^\ell)$ be the space of $1$-jets of 
mappings of $N$ into $\R^\ell$. 
Then, note that $J^1(N,\R^\ell)$ is a $C^{r-1}$ manifold. 
For a given $C^r$ mapping $g:N\to \R^\ell$ $(r\geq 2)$, 
the mapping $j^1 g:N\to J^1(N,\R^\ell)$ 
is defined by $q \mapsto j^1 g(q)$. 
Then, notice that the mapping $j^1 g:N\to J^1(N,\R^\ell)$ is of class $C^{r-1}$. 
For details on the space $J^1(N,\R^\ell)$ or 
the mapping $j^1 g:N\to J^r(N,\R^\ell)$, see for example, \cite{GG}. 

Now, let $\{(U_\lambda ,\varphi _\lambda )\}_{\lambda \in \Lambda}$ be a coordinate neighborhood system of $N$. 
Let $\Pi :J^1(N,\mathbb{R}^\ell)$$\to N\times \mathbb{R}^\ell$ be the natural projection defined by $\Pi(j^1g(q))=(q,g(q))$. 
Let $\Phi _\lambda :\Pi^{-1}(U_\lambda \times \mathbb{R}^\ell)
\to \varphi _\lambda (U_\lambda)\times \mathbb{R}^\ell \times J^1(n,\ell)$ be 
the homeomorphism defined by 
\begin{eqnarray*} 
\Phi _\lambda \left(j^1g(q)\right)=\left(\varphi _\lambda (q),g(q),j^1(\psi_{_\lambda} \circ g\circ \varphi _\lambda ^{-1}\circ \widetilde{\varphi} _\lambda)(0)\right), 
\end{eqnarray*} 
where 
$J^1(n, \ell)=\{ j^1g(0) \mid g : (\mathbb{R}^n, 0) \to  (\mathbb{R}^\ell, 0) \}$ 
and 
$\widetilde{\varphi} _\lambda : \mathbb{R}^n\to \mathbb{R}^n$ 
(resp., $\psi_{\lambda} : \mathbb{R}^m \to \mathbb{R}^m $) is the 
translation given by 
$\widetilde{\varphi} _\lambda(0)=\varphi _\lambda (q)$ 
(resp., $\psi_{\lambda}(g(q))=0$). 
Then, $\{(\Pi^{-1}(U_\lambda \times \mathbb{R}^\ell), 
\Phi _\lambda )\}_{\lambda \in \Lambda}$ 
is a coordinate neighborhood system of $J^1(N,\mathbb{R}^\ell)$. 
Set 
\begin{eqnarray*} 
\Sigma ^k=\left\{j^1g(0)\in J^1(n,\ell)\mid {\rm corank\ }Jg(0)=k\right\}, 
\end{eqnarray*}
where ${\rm corank\ }Jg(0)={\rm min} \{n,\ell \}-{\rm rank\ } Jg(0)$ 
and $k=1,2,\ldots ,{\rm min}\{n, \ell\}$. 
Set 
\begin{eqnarray*}
\Sigma^k(N,\mathbb{R}^\ell)=\bigcup_{\lambda \in \Lambda}\Phi ^{-1}_\lambda \left(\varphi _\lambda (U_\lambda )\times \mathbb{R}^\ell \times \Sigma ^k\right). 
\end{eqnarray*}
Then, the set $\Sigma ^k(N,\mathbb{R}^\ell)$ is a submanifold of 
$J^1(N,\mathbb{R}^\ell)$ satisfying  
\begin{eqnarray*}
{\rm codim}\ \Sigma ^k(N,\mathbb{R}^\ell)&=&{\rm dim}\ J^1(N,\mathbb{R}^\ell)- 
{\rm dim}\ \Sigma ^k(N,\mathbb{R}^\ell) \\
&=&(n-v+k)(\ell-v+k), 
\end{eqnarray*}
where $v={\rm min}\{n, \ell\}$. 
(For details on $\Sigma^k$ and $\Sigma^k(N,\mathbb{R}^\ell)$, see for instance 
\cite{GG}, pp. 60--61). 

Then, the first main theorem in this paper is the following. 
\begin{theorem}\label{main} 
Let $f$ be a $C^r$ immersion 
of $N$ into an open subset $U$ of $\mathbb{R}^m$, 
where $N$ is a $C^r$ manifold of dimension $n$. 
Let $F:U\to \mathbb{R}^\ell$ be a $C^r$ mapping and 
$k$ be a positive integer satisfying $1\leq k \leq \min \{n,\ell \}$. 
If 
\[
r>\max \{ \dim N -\codim \Sigma ^k(N,\mathbb{R}^\ell), 0\}+1, 
\]
then there exists a subset $\Sigma$ with Lebesgue measure zero 
of $\mathcal{L}(\mathbb{R}^{m},\mathbb{R}^{\ell})$ 
such that for any $\pi \in \mathcal{L}(\mathbb{R}^{m},\mathbb{R}^{\ell})-\Sigma$, 
the mapping $j^1(F_\pi \circ f):
N\to J^1(N,\mathbb{R}^\ell)$ is transverse to the 
submanifold $\Sigma^k(N,\mathbb{R}^\ell)$.
\end{theorem}

Now, in order to state the second main theorem (\cref{main2}), 
we will prepare some definitions. 
Let $N$ be a $C^r$ manifold $(r\geq1)$. 
Set 
\[
N^{(s)}=\{(q_1,\<q_2,\>\ldots ,q_s)\in N^s\mid q_i\not=q_j \ (i\not=j)\}.
\] 
Note that $N^{(s)}$ is an open submanifold of $N^s$. 
For any mapping $g:N\to \R^\ell$, 
let $g^{(s)}:N^{(s)}\to (\R^\ell)^s$ be the mapping given by 
\begin{eqnarray*}
g^{(s)}(q_1,\<q_2,\>\ldots ,q_s)=(g(q_1),\<g(q_2),\>\ldots ,g(q_s)). 
\end{eqnarray*}
Set 
$\Delta_s=\{(y,\ldots ,y) \in (\R^\ell)^s \mid y\in \R^\ell \}$. 
Then, $\Delta_s$ is a submanifold of $(\R^\ell)^s$ satisfying 
\begin{eqnarray*}
\codim \Delta_s =\dim (\R^\ell)^s- 
\dim \Delta_s  =\ell (s-1).
\end{eqnarray*}

\begin{definition}\label{normal crossings}
{\rm
Let $g$ be a $C^1$ mapping of $N$ into $\R^\ell$, where 
$N$ is a $C^r$ manifold $(r\geq 1)$.   
Then, $g$ is called a {\it mapping with normal crossings} if  
for any positive integer $s$ $(s\geq 2)$, 
the mapping $g^{(s)}:N^{(s)}\to (\R^\ell)^s$ is 
transverse to $\Delta_s$. 
}
\end{definition}
As in \cite{CG}, for any injection $f:N\to \mathbb{R}^m$, set 
\begin{eqnarray*}
s_f={\rm max}\left\{s\ \middle| \ \forall (q_1,\<q_2,\>\ldots ,q_s)\in N^{(s)}, {\rm dim\ }\sum_{i=2}^s \mathbb{R}\overrightarrow{f(q_1)f(q_i)}=s-1 \right\}. 
\end{eqnarray*}
Since the mapping $f$ is an injection, we have $2 \leq s_f$. 
Since $f(q_1),f(q_2),\ldots ,f(q_{s_f})$ are points of $\mathbb{R}^m$, 
it follows that $s_f\leq m+1$. 
Hence, we get  
\begin{eqnarray*}
2\leq s_f \leq m+1. 
\end{eqnarray*}
Moreover, in the following, for a set $X$, 
we denote the number of its elements (or its cardinality) by $|X|$. 
Then, the second main theorem in this paper is the following. 
\begin{theorem}\label{main2}
Let $f$ be a $C^r$ injection of $N$ into an open subset 
$U$ of $\mathbb{R}^m$, 
where $N$ is a $C^r$ manifold of dimension $n$. 
Let $F:U\to \mathbb{R}^\ell$ be a $C^r$ mapping. 
If 
\[
r>\max \{ s_0,0\},
\]
then there exists a subset $\Sigma$ of $\mathcal{L}(\mathbb{R}^m, \mathbb{R}^\ell)$ with Lebesgue measure zero such that 
for any $\pi \in \mathcal{L}(\mathbb{R}^m, \mathbb{R}^\ell)-\Sigma $, and  
for any $s$ $(2\leq s \leq s_f)$, 
the $C^r$ mapping $(F_\pi \circ f)^{(s)}:N^{(s)}\to (\R^\ell)^s$ is transverse to the 
submanifold $\Delta_s$, 
where 
\[
s_0=\max \{ s(n-\ell)+\ell \ |\ 2 \leq s \leq s_f\}. 
\]
Moreover, if the mapping $F_\pi$ 
satisfies that $| F_\pi^{-1}(y) | \leq s_f$ for any $y\in\mathbb{R}^\ell$, 
then $F_\pi \circ f:N\to \mathbb{R}^\ell$ is a  $C^r$ mapping with normal crossings.    
\end{theorem}
\begin{remark}
{\rm 
\begin{enumerate}
\item 
There is an advantage that the domain of the mapping $F$ 
is an arbitrary open set. 
Suppose that $U=\mathbb{R}$. 
Let $F:\mathbb{R}\to \mathbb{R}$ be the function defined by $x\mapsto |x|$. 
Since $F$ is not differentiable at $x=0$, we cannot apply \cref{main,main2} to 
$F:\mathbb{R}\to \mathbb{R}$. 

On the other hand, if $U=\mathbb{R}-\{0\}$, then 
\cref{main,main2} 
can be applied to the restriction $F|_U$.
\item 
As in \cite{CG}, 
there is a case of $s_f =3$ as follows. 
If $n+1 \leq m$, $N=S^n$ and 
$f:S^n\to \mathbb{R}^m$ is the inclusion $f(x)=(x,0,\ldots ,0)$, 
then we get $s_f=3$. 
Indeed, suppose that there exists a point $(q_1, q_2, q_3)\in (S^n)^{(3)}$ 
satisfying ${\rm dim\ }\sum_{i=2}^3 \mathbb{R}\overrightarrow{f(q_1)f(q_i)}=1$. 
Then, since the number of the intersections of $f(S^n)$ and 
a straight line of $\mathbb{R}^m$ is at most two, 
this contradicts the assumption. 
Thus, 
we have $s_f\geq3$. 
From $S^1\times \{0\} \subset f(S^n)$, 
we get $s_f<4$, where $0=\underbrace{(0,\ldots ,0)}_{(m-2)
\text{-tuple}}$. 
Therefore, it follows that $s_f=3$.   
\item 

The essential idea for the proofs of \cref{main,main2}   
is to apply \cref{abra}, and it 
is similar to the idea of the proofs of \cite[Theorems~1~and~2]{CG}. 
Note that 
in the special case $r=\infty$, from some results in \cite{CG}, the 
results in this paper 
(\cref{main,main2} in this section and \cref{Morse function,immersion,corank,normal,injective,injective immersion,embedding} 
in \cref{section 5}) 
can be obtained. 
\end{enumerate}
}
\end{remark}

The following well known result is important for the proofs of \cref{main,main2}. 
In \cite{GG}, the proof of \cref{abra} in the case $r=\infty$ is given. 
Hence, for the sake of readers' convenience, 
the proof of \cref{abra} is given in \cref{app} as an appendix. 
\begin{lemma}[\cite{GG}]\label{abra}
Let $N$, $A$, $P$ be $ C^r$ manifolds, 
$Z$ be a $C^r$ submanifold of $P$ and $\Gamma:N\times A\to P$ be a 
$C^r$ mapping. 
If 
\[
r>\max \{ \dim N-\codim Z,0\},
\]
and $\Gamma$ is transverse to $Z$, then there exists a subset $\Sigma$ of $A$ 
with Lebesgue measure zero such that for any $a \in A-\Sigma $, 
the $C^r$ mapping 
$\Gamma_a:N\to P$ is transverse to $Z$, 
where $\codim Z=\dim P-\dim Z$ and $\Gamma_a(q)=\Gamma(q,a)$. 
\end{lemma}

\section{Proof of Theorem \ref{main}}\label{section 3}
In this proof, for a positive integer $\wt{n}$, 
we denote the $\wt{n}\times \wt{n}$ unit matrix by $E_{\wt{n}}$. 
Let $(\alpha_{ij})_{1\leq i \leq \ell, 1\leq j \leq m}$ be a representing matrix of 
a linear mapping $\pi:\mathbb{R}^m\to \mathbb{R}^\ell$. 
Set $F_{\alpha}=F_{\pi}$. Then, we have
\[
F_{\alpha}(x)=
\biggl(F_{1}(x)+\sum_{j=1}^{m}\alpha_{1j}x_j, 
F_{2}(x)+\sum_{j=1}^{m}\alpha_{2j}x_j, 
\ldots,
F_{\ell}(x)+\sum_{j=1}^{m}\alpha_{\ell j}x_j\biggr), 
\eqno (3.1)
\]
where $F=(F_1, F_2,\ldots ,F_\ell)$, 
 $\alpha=(\alpha_{11},\<\alpha_{12},\ldots, 
 \alpha_{1m},\ldots ,
 \alpha_{\ell 1},\alpha_{\ell 2},\ldots ,
\alpha_{\ell m})\in (\mathbb{R}^m)^\ell$ and $x=(x_1,x_2,\ldots ,x_m)$. 
For a given $C^r$ immersion $f:N\to U$, 
the $C^r$ mapping $F_{\alpha}\circ f:N\to \mathbb{R}^\ell$ is 
given as follows: 
\[
F_{\alpha}\circ f=
\biggl(F_{1}\circ f+\sum_{j=1}^{m}\alpha_{1j}f_j, 
F_{2}\circ f+\sum_{j=1}^{m}\alpha_{2j}f_j, 
\ldots,
F_{\ell}\circ f+\sum_{j=1}^{m}\alpha_{\ell j}f_j\biggr),
\eqno (3.2)
\]
where $f=(f_1,f_2,\ldots ,f_m)$.  
Since we have the natural identification 
$\mathcal{L}(\mathbb{R}^{m},\mathbb{R}^{\ell})=(\mathbb{R}^{m})^\ell$, 
for the proof, 
it is sufficient to show that 
there exists a subset $\Sigma$ 
with Lebesgue measure zero of $(\mathbb{R}^m)^{\ell}$ 
such that for any $\alpha \in (\mathbb{R}^m)^{\ell}-\Sigma$, 
the mapping $j^{1}(F_{\alpha}\circ f):
N\to J^{1}(N,\mathbb{R}^\ell)$ is transverse to $\Sigma^k(N,\mathbb{R}^\ell)$. 

Now, let $\Gamma:N\times (\mathbb{R}^m)^\ell \to J^1(N,\mathbb{R}^\ell)$ 
be the $C^{r-1}$ mapping defined by 
\[
\Gamma(q,\alpha)=j^1(F_\alpha \circ f)(q).
\]
Note that $r-1>\max\{\dim N-\codim \Sigma^k(N,\R^\ell),0\}$. 
Thus, 
if $\Gamma$ is transverse to 
$\Sigma^k(N,\mathbb{R}^\ell)$, 
then from \cref{abra}, 
there exists a subset $\Sigma$ of $(\mathbb{R}^m)^\ell$ with 
Lebesgue measure zero such that 
for any $\alpha \in (\mathbb{R}^m)^\ell-\Sigma $, 
the $C^{r-1}$ mapping $\Gamma_\alpha :N\to J^1(N,\mathbb{R}^\ell)$ 
($\Gamma_\alpha=j^1(F_\alpha \circ f)$) 
is transverse to 
$\Sigma^k(N,\mathbb{R}^\ell)$. 
Therefore, for the proof, 
it is sufficient to show that if $\Gamma(\widetilde{q},\widetilde{\alpha})
\in \Sigma^k(N,\mathbb{R}^\ell)$, then the following holds:
\[
d\Gamma_{(\widetilde{q},\widetilde{\alpha})}(T_{(\widetilde{q},\widetilde{\alpha})}
(N\times (\mathbb{R}^m)^\ell))
+
T_{\Gamma(\widetilde{q},\widetilde{\alpha})}\Sigma^k(N,\mathbb{R}^\ell)
=
T_{\Gamma(\widetilde{q},\widetilde{\alpha})}J^1(N,\mathbb{R}^\ell). \eqno (3.3)
\]
As in \cref{section 2}, 
let $\{(U_\lambda ,\varphi _\lambda )\}_{\lambda \in \Lambda}$ 
(resp., $\{(\Pi^{-1}(U_\lambda \times \mathbb{R}^\ell), 
\Phi _\lambda )\}_{\lambda \in \Lambda}$) 
be a coordinate neighborhood system of $N$ (resp., $J^1(N,\mathbb{R}^\ell)$). 
There exists a coordinate neighborhood 
$\left(U_{\widetilde{\lambda}}\times (\mathbb{R}^m)^\ell, \varphi_{\widetilde{\lambda}}\times id \right)$ 
containing the point $(\widetilde{q},\widetilde{\alpha})$ of 
$N\times (\mathbb{R}^m)^\ell$, 
where $id$ is the identity mapping of 
$(\mathbb{R}^m)^\ell$ into $(\mathbb{R}^m)^\ell$, 
and the mapping $\varphi_{\widetilde{\lambda}}\times id : 
U_{\widetilde{\lambda}}\times (\mathbb{R}^m)^\ell \to 
\varphi_{\widetilde{\lambda}}(U_{\widetilde{\lambda}})\times (\mathbb{R}^m)^\ell$ ($\subset \mathbb{R}^n\times (\mathbb{R}^m)^\ell$) 
is given by 
$\left(\varphi_{\widetilde{\lambda}}\times id\right)(q,\alpha)=
\left(\varphi_{\widetilde{\lambda}}(q), id(\alpha)\right)$.
There exists a coordinate neighborhood 
$\left(\Pi^{-1}(U_{\widetilde{\lambda}} \times \mathbb{R}^\ell), 
\Phi _{\widetilde{\lambda}} \right)$ 
containing the point $\Gamma (\widetilde{q},\widetilde{\alpha})$ of 
$J^1(N,\mathbb{R}^\ell)$. 
Let $t=(t_1,\<t_2,\>\ldots ,t_n)\in \mathbb{R}^n$ be a local coordinate on 
$\varphi_{\widetilde{\lambda}}(U_{\widetilde{\lambda}})$ 
containing 
$\varphi_{\widetilde{\lambda}}(\widetilde{q})$. 
Then, the mapping $\Gamma$ is locally given by the following:
{\small 
\begin{eqnarray*}
&{}&(\Phi _{\widetilde{\lambda}} \circ \Gamma \circ 
(\varphi_{\widetilde{\lambda}}\times id )^{-1})(t,\alpha)\\
&=&(\Phi _{\widetilde{\lambda}}\circ j^1(F_\alpha \circ f)\circ \varphi_{\widetilde{\lambda}}^{-1})(t)\\
&=&\left(t, (F_\alpha \circ f \circ \varphi_{\widetilde{\lambda}}^{-1})(t), \right.  \\
&& \frac{\partial (F_{\alpha, 1}\circ f\circ \varphi_{\widetilde{\lambda}}^{-1})}
{\partial t_1}(t), 
\<
\frac{\partial (F_{\alpha, 1}\circ f\circ \varphi_{\widetilde{\lambda}}^{-1})}
{\partial t_2}(t), 
\>
\ldots ,
\frac{\partial (F_{\alpha, 1}\circ f\circ \varphi_{\widetilde{\lambda}}^{-1})}
{\partial t_n}(t),
\\ 
&& \< 
\frac{\partial (F_{\alpha, 2}\circ f\circ \varphi_{\widetilde{\lambda}}^{-1})}
{\partial t_1}(t), 
\frac{\partial (F_{\alpha, 2}\circ f\circ \varphi_{\widetilde{\lambda}}^{-1})}
{\partial t_2}(t), 
\ldots ,
\frac{\partial (F_{\alpha, 2}\circ f\circ \varphi_{\widetilde{\lambda}}^{-1})}
{\partial t_n}(t), 
\> 
\\ 
\\
&&\hspace{120pt}\cdots \cdots \cdots , \\ 
\\
&& \left. \frac{\partial (F_{\alpha, \ell}\circ f\circ \varphi_{\widetilde{\lambda}}^{-1})}
{\partial t_1}(t)
,
\<
\frac{\partial (F_{\alpha, \ell}\circ f\circ \varphi_{\widetilde{\lambda}}^{-1})}
{\partial t_2}(t), 
\>
\ldots ,
\frac{\partial (F_{\alpha, \ell}\circ f\circ \varphi_{\widetilde{\lambda}}^{-1})}
{\partial t_n}(t)
\right)\\
&=&\left(t, (F_\alpha \circ f \circ \varphi_{\widetilde{\lambda}}^{-1})(t), \right. \\
&&\frac{\partial F_1\circ \widetilde{f}}{\partial t_1}(t)+
\sum_{j=1}^m \alpha_{1j}\frac{\partial \widetilde{f}_j}{\partial t_1}(t),
\<
\frac{\partial F_1\circ \widetilde{f}}{\partial t_2}(t)+
\sum_{j=1}^m \alpha_{1j}\frac{\partial \widetilde{f}_j}{\partial t_2}(t),
\>
\ldots ,
\frac{\partial F_1\circ \widetilde{f}}{\partial t_n}(t)+
\sum_{j=1}^m \alpha_{1j}\frac{\partial \widetilde{f}_j}{\partial t_n}(t),
\\ 
&& 
\<
\frac{\partial F_2\circ \widetilde{f}}{\partial t_1}(t)+
\sum_{j=1}^m \alpha_{2j}\frac{\partial \widetilde{f}_j}{\partial t_1}(t),
\frac{\partial F_2\circ \widetilde{f}}{\partial t_2}(t)+
\sum_{j=1}^m \alpha_{2j}\frac{\partial \widetilde{f}_j}{\partial t_2}(t),
\ldots ,
\frac{\partial F_2\circ \widetilde{f}}{\partial t_n}(t)+
\sum_{j=1}^m \alpha_{2j}\frac{\partial \widetilde{f}_j}{\partial t_n}(t),
\> 
\\
\\
&&\hspace{160pt}\cdots \cdots \cdots , 
\\
\\
&& \left. 
\frac{\partial F_\ell \circ \widetilde{f}}{\partial t_1}(t)+\sum_{j=1}^m \alpha_{\ell j}\frac{\partial \widetilde{f}_j}{\partial t_1}(t),
\<
\frac{\partial F_\ell \circ \widetilde{f}}{\partial t_2}(t)+\sum_{j=1}^m \alpha_{\ell j}\frac{\partial \widetilde{f}_j}{\partial t_2}(t),
\>
\ldots ,
\frac{\partial F_\ell \circ \widetilde{f}}{\partial t_n}(t)+
\sum_{j=1}^m \alpha_{\ell j}\frac{\partial \widetilde{f}_j}{\partial t_n}(t)\right), 
\end{eqnarray*} 
}where $F_\alpha=(F_{\alpha, 1},\<F_{\alpha, 2},\> \ldots , F_{ \alpha, \ell})$ and 
$\widetilde{f}=(\widetilde{f}_1, \<\widetilde{f}_2,\> \ldots ,\widetilde{f}_m)
=(f_1\circ \varphi_{\widetilde{\lambda}}^{-1},
\<
f_2\circ \varphi_{\widetilde{\lambda}}^{-1},
\>
\ldots , f_m\circ
\varphi_{\widetilde{\lambda}}^{-1})=f\circ \varphi_{\widetilde{\lambda}}^{-1}$. 
The Jacobian matrix of $\Gamma$ at 
$(\widetilde{q},\widetilde{\alpha})$ is the following:\\
\begin{eqnarray*}
J\Gamma_{(\widetilde{q},\widetilde{\alpha})}=
\left(
\begin{array}{@{\,}c@{\,\,}|@{\,\,}c@{\,\,\,}c@{\,\,\,}c@{\,\,\,}c@{\,\,\,}c}
E_n       &    0    &         \cdots          &   \cdots        &0    \\
\hline                &   \ast      &    \cdots    &  \cdots      & \ast  \\ 
      &   {}^t\!(Jf_{\widetilde{q}}) & & \bigzerol &\\
 \ast    &              & {}^t\!(Jf_{\widetilde{q}})    &       &  \\
    &    & \bigzerol &\ddots &      \\
      &    &  &  &   {}^t\!(Jf_{\widetilde{q}}) \\
\end{array}
\right)_{(t, \alpha)=(\varphi_{\widetilde{\lambda}}(\widetilde{q}),\widetilde{\alpha})}, 
\end{eqnarray*}
where 
$Jf_{\widetilde{q}}$ is the Jacobian matrix of $f$ at $\widetilde{q}$. 
Notice that ${}^t\!(Jf_{\widetilde{q}})$ is the transpose of 
$Jf_{\widetilde{q}}$ and that there are $\ell$ copies of ${}^t\!(Jf_{\widetilde{q}})$ 
in the above description of $J\Gamma_{(\widetilde{q},\widetilde{\alpha})}$. 
Since $\Sigma^k(N,\mathbb{R}^\ell)$ is a subfiber-bundle of 
$J^1(N,\mathbb{R}^\ell)$ with the fiber $\Sigma^k$, in order to show (3.3), 
it is sufficient to prove that the matrix $M_1$ given below 
has rank $n+\ell+n\ell$: 
\begin{eqnarray*}
M_1=
\left(
\begin{array}{@{\,}c@{\,\,}|@{\,\,}c@{\,\,\,}c@{\,\,\,}c@{\,\,\,}c@{\,\,\,}c}
E_{n+\ell}       &    \ast    &         \cdots          &   \cdots        &\ast   \\
\hline     \\[-3.3mm]   &   {}^t\!(Jf_{\widetilde{q}})  & & \bigzerol &\\
\bigzerol  &              & {}^t\!(Jf_{\widetilde{q}})      &       &  \\
    &    & \bigzerol &\ddots &      \\
      &    &  &  &   {}^t\!(Jf_{\widetilde{q}})  \\
\end{array}
\right)_{(t, \alpha)=(\varphi_{\widetilde{\lambda}}(\widetilde{q}),\widetilde{\alpha})}.
\end{eqnarray*}
Notice that there are $\ell$ copies of ${}^t\!(Jf_{\widetilde{q}})$ 
in the above description of $M_1$. 
Note that for any $i$ $(1\leq i \leq m\ell)$, the $(n+\ell+i)$-th 
column vector of $M_1$ coincides with 
the $(n+i)$-th column vector of $J\Gamma_{(\widetilde{q},\widetilde{\alpha})}$. 
Since $f$ is an immersion $(n\leq m)$, the rank of 
$M_1$ 
is equal to $n+\ell+n\ell$. 
Therefore, we get (3.3). 
\hfill\qed
\section{Proof of \cref{main2}}\label{section 4}
As in the proof of \cref{main}, 
set $F_{\alpha}=F_{\pi}$, where $F_{\alpha}$ is given by (3.1) 
in \cref{section 3}. 
For a given $C^r$ injection $f:N\to U$, 
the $C^r$ mapping $F_{\alpha}\circ f:N\to \mathbb{R}^\ell$ is 
given by the same expression as (3.2).  
Since we have the natural identification 
$\mathcal{L}(\mathbb{R}^{m},\mathbb{R}^{\ell})=(\mathbb{R}^{m})^\ell$, 
in order to prove that 
there exists a subset $\Sigma$ of $\mathcal{L}(\mathbb{R}^m, \mathbb{R}^\ell)$ with Lebesgue measure zero such that 
for any $\pi \in \mathcal{L}(\mathbb{R}^m, \mathbb{R}^\ell)-\Sigma $, 
\<and \> 
for any $s$ $(2\leq s \leq s_f)$, 
the $C^r$ mapping $(F_\pi \circ f)^{(s)}:N^{(s)}\to (\R^\ell)^s$ is transverse to 
$\Delta_s$,  
it is sufficient to prove that 
there exists a subset $\Sigma$ of $(\mathbb{R}^m)^{\ell}$ 
with Lebesgue measure zero 
such that for any $\alpha \in (\mathbb{R}^m)^{\ell}-\Sigma$, 
\<and \> 
for any $s$ $(2\leq s \leq s_f)$, 
the $C^r$ mapping $(F_{\alpha}\circ f)^{(s)}:
N^{(s)}\to (\mathbb{R}^\ell)^s$ is transverse to the submanifold $\Delta_s$.  

Now, let $s$ be a positive integer satisfying $2\leq s \leq s_f$.  
Let $\Gamma : N^{(s)}\times (\mathbb{R}^m)^\ell \to (\mathbb{R}^{\ell})^s$ be 
the $C^r$ mapping given by 
\begin{eqnarray*}
\Gamma(q_1,\<q_2,\>\ldots ,q_s,\alpha)=
\left( (F_\alpha \circ f)(q_1), \< (F_\alpha \circ f)(q_2),\>\ldots, 
(F_\alpha \circ f)(q_s)\right).
\end{eqnarray*}
Note that from $r>\max\{s_0,0\}$, 
we have 
\begin{eqnarray*}
r&>&\max\{s(n-\ell)+\ell,0\}
\\
&=&\max\{\dim N^{(s)}-\codim \Delta_s,0\}
\end{eqnarray*}
for any positive integer $s$ $(2\leq s \leq s_f)$. 
Thus, if for any positive integer $s$ $(2\leq s \leq s_f)$, 
the mapping $\Gamma$ is transverse to $\Delta_s$, then from \cref{abra}, 
for any positive integer $s$ $(2\leq s \leq s_f)$, 
 there exists a subset $\Sigma_s$ of $(\mathbb{R}^m)^\ell$ 
 with Lebesgue measure zero 
such that for any $\alpha  \in (\mathbb{R}^m)^{\ell}-\Sigma_s$,  
the mapping $\Gamma_\alpha : 
N^{(s)}\to (\mathbb{R}^\ell)^s$ $(\Gamma_\alpha=(F_\alpha\circ f)^{(s)})$ 
is transverse to $\Delta_s$. 
Then, $\Sigma=\bigcup_{s=2}^{s_f}\Sigma_s$ is a subset of  $(\mathbb{R}^m)^\ell$ 
with Lebesgue measure zero. 
Thus, 
for any $\alpha \in (\mathbb{R}^m)^{\ell}-\Sigma $, and  
for any $s$ $(2\leq s \leq s_f)$, 
the $C^r$ mapping $\Gamma_\alpha : 
N^{(s)}\to (\mathbb{R}^\ell)^s$ $(\Gamma_\alpha=(F_\alpha\circ f)^{(s)})$ 
is transverse to $\Delta_s$. 

Therefore, for this proof, 
it is sufficient to prove that 
for any positive integer $s$ $(2\leq s \leq s_f)$,   
if $\Gamma(\widetilde{q}, \widetilde{\alpha})\in \Delta_s $ 
$(\widetilde{q}=(\widetilde{q}_1,
\widetilde{q}_2,
\ldots ,\widetilde{q}_s))$, 
then the following holds: 
\[
d\Gamma_{(\widetilde{q}, \widetilde{\alpha})}(T_{(\widetilde{q}, \widetilde{\alpha})}(N^{(s)}\times (\mathbb{R}^m)^\ell))+
T_{\Gamma(\widetilde{q}, \widetilde{\alpha})}\Delta_s
=T_{\Gamma(\widetilde{q}, \widetilde{\alpha})}(\mathbb{R}^\ell)^s. \eqno (4.1)
\]
Let $\{(U_\lambda ,\varphi _\lambda )\}_{\lambda \in \Lambda}$ be a coordinate  neighborhood system of $N$. 
There exists a coordinate neighborhood 
$(U_{\widetilde{\lambda}_1}\times \<U_{\widetilde{\lambda}_2}\times \>
\cdots \times U_{\widetilde{\lambda}_s} \times (\mathbb{R}^m)^\ell, 
\varphi_{\widetilde{\lambda}_1}\times 
\<\varphi_{\widetilde{\lambda}_2}\times 
\>
\cdots \times \varphi_{\widetilde{\lambda}_s}\times id )$ 
containing $(\widetilde{q}, \widetilde{\alpha})$ 
of $N^{(s)}\times (\mathbb{R}^m)^\ell$, 
where $id:(\mathbb{R}^m)^\ell\to (\mathbb{R}^m)^\ell$ is the identity mapping, 
and 
$\varphi_{\widetilde{\lambda}_1}\times 
\<
\varphi_{\widetilde{\lambda}_2}\times 
\>
\cdots \times \varphi_{\widetilde{\lambda}_s}\times id  : 
U_{\widetilde{\lambda}_1}\times 
\<
U_{\widetilde{\lambda}_2}\times 
\>
\cdots \times U_{\widetilde{\lambda}_s} \times (\mathbb{R}^m)^\ell 
 \to (\mathbb{R}^n)^s\times (\mathbb{R}^m)^\ell$ is defined by 
$(\varphi_{\widetilde{\lambda}_1}\times 
\<
\varphi_{\widetilde{\lambda}_2}\times 
\>
\cdots \times \varphi_{\widetilde{\lambda}_s}\times id )(q_1,\<q_2,\>\ldots ,q_s,\alpha)=
(\varphi_{\widetilde{\lambda}_1}(q_1),  
\<
\varphi_{\widetilde{\lambda}_2}(q_2),  
\>
\ldots ,
\varphi_{\widetilde{\lambda}_s}(q_s), id(\alpha ))$.
Let $t_i=(t_{i1},t_{i2},\ldots ,t_{in})$ be a local coordinate 
around $\varphi_{\widetilde{\lambda}_i}(\widetilde{q}_i)$ $(1\leq i \leq s)$.  
Then, $\Gamma$ is locally given by the following:
\begin{eqnarray*}
&&\Gamma \circ \left(\varphi_{\widetilde{\lambda}_1}\times 
\<
\varphi_{\widetilde{\lambda}_2}\times 
\>
\cdots \times \varphi_{\widetilde{\lambda}_s}\times id\right)^{-1}
(t_1,\<t_2,\>\ldots ,t_s,\alpha)
\\
&=&\left( (F_\alpha \circ f\circ \varphi_{\widetilde{\lambda}_1}^{-1})(t_1), 
\<
(F_\alpha \circ f\circ \varphi_{\widetilde{\lambda}_2}^{-1})(t_2), 
\>
\ldots ,
(F_\alpha \circ f\circ \varphi_{\widetilde{\lambda}_s}^{-1})(t_s) 
\right)
\\
&=&\left(
F_1\circ \widetilde{f}(t_1)+\sum_{j=1}^m\alpha_{1j}\widetilde{f}_j(t_1), 
\<
F_2\circ \widetilde{f}(t_1)+\sum_{j=1}^m\alpha_{2j}\widetilde{f}_j(t_1), 
\>
\ldots ,
F_\ell \circ \widetilde{f}(t_1)+\sum_{j=1}^m\alpha_{\ell j}\widetilde{f}_j(t_1), \right.\\
&& 
\<
F_1\circ \widetilde{f}(t_2)+\sum_{j=1}^m\alpha_{1j}\widetilde{f}_j(t_2), 
F_2\circ \widetilde{f}(t_2)+\sum_{j=1}^m\alpha_{2j}\widetilde{f}_j(t_2), 
\ldots ,
F_\ell \circ \widetilde{f}(t_2)+\sum_{j=1}^m\alpha_{\ell j}\widetilde{f}_j(t_2), 
\>  
\\
\\
&&\hspace{150pt}\cdots \cdots \cdots , 
\\
\\
&& \left. 
F_1\circ \widetilde{f}(t_s)+\sum_{j=1}^m\alpha_{1j}\widetilde{f}_j(t_s),
\<
F_2\circ \widetilde{f}(t_s)+\sum_{j=1}^m\alpha_{2j}\widetilde{f}_j(t_s),
\>
\ldots ,
F_\ell \circ \widetilde{f}(t_s)+\sum_{j=1}^m\alpha_{\ell j}\widetilde{f}_j(t_s)
\right), 
\end{eqnarray*}  
where
$\widetilde{f}(t_i)=(\widetilde{f}_1(t_i),
\<
\widetilde{f}_2(t_i),
\>
\ldots ,\widetilde{f}_m(t_i))
=(f_1\circ \varphi_{\widetilde{\lambda}_i}^{-1}(t_i),
\<
f_2\circ \varphi_{\widetilde{\lambda}_i}^{-1}(t_i),
\>
\ldots ,
f_m\circ \varphi_{\widetilde{\lambda}_i}^{-1}(t_i))$ $(1\leq i\leq s)$. 
For simplicity, set $t=(t_1,\<t_2,\>\ldots ,t_s)$ and 
$z=(\varphi_{\widetilde{\lambda}_1}\times 
\<
\varphi_{\widetilde{\lambda}_2}\times 
\>
\cdots \times \varphi_{\widetilde{\lambda}_s})
(\widetilde{q}_1,\<\widetilde{q}_2,\>\ldots ,\widetilde{q}_s)$. 

The Jacobian matrix of $\Gamma$ at 
$(\widetilde{q}, \widetilde{\alpha})$ is the following: 
\begin{eqnarray*}
J\Gamma_{(\widetilde{q}, \widetilde{\alpha})}=
\left(
\begin{array}{@{\,\,\,}c@{\,\,\,\,\,}|@{\,\,\,\,\,}c@{\,\,\,}}
\ast & B(t_1) \\
\ast & B(t_2) \\
 \vdots & \vdots    \\ 
\ast & B(t_s) \\
\end{array}
\right)_{(t,\alpha)=(z, \widetilde{\alpha})},
\end{eqnarray*}
where 
\begin{eqnarray*}
B(t_i)=
\left. 
\left(
\begin{array}{ccccccc}
{\bf b}_{}(t_i)& & &\bigzerol \\
&{\bf b}_{}(t_i)&& \\
\bigzerol &   &  \ddots   &    \\ 
&&   &         &{\bf b}_{}(t_i)
\end{array}
\right)
\right\}
\,\text{{\rm $\ell$ rows}}
\end{eqnarray*}
and 
${\bf b}_{ }(t_i)=(
\widetilde{f}_1(t_i),
\widetilde{f}_2(t_i),
\ldots ,\widetilde{f}_m(t_i))$. 
By the construction of 
$T_{\Gamma(\widetilde{q}, \widetilde{\alpha})}\Delta_s$, 
in order to prove (4.1), it is sufficient to 
prove that the rank of the following matrix $M_2$ is equal to $\ell s$: 
\begin{eqnarray*}
M_2=
\left(
\begin{array}{@{\,\,\,}c@{\,\,\,\,\,}|@{\,\,\,\,\,}c@{\,\,\,}}
E_\ell & B(t_1) \\
E_\ell & B(t_2) \\
 \vdots & \vdots    \\ 
E_\ell & B(t_s) \\
\end{array}
\right)_{t=z}.
\end{eqnarray*}
There exists an 
$\ell s\times \ell s$ regular matrix $Q_1$ 
satisfying 
\begin{eqnarray*}
Q_1M_2=
\left(
\begin{array}{@{\,\,\,}c@{\,\,\,\,\,}|@{\,\,\,\,\,}c@{\,\,\,}}
E_\ell & B(t_1) \\
0 & B(t_2)- B(t_1) \\
 \vdots & \vdots    \\ 
0 & B(t_s)- B(t_1) \\
\end{array}
\right)_{t=z}.
\end{eqnarray*}
There exists an 
$(\ell+m\ell)\times (\ell+m\ell)$ regular matrix $Q_2$ 
satisfying  
\begin{eqnarray*}
Q_1M_2Q_2&=&
\left(
\begin{array}{@{\,\,\,}c@{\,\,\,\,}|@{\,\,\,\,}c@{\,\,\,}}
E_\ell & 0 \\
0 & B(t_2)- B(t_1) \\
 \vdots & \vdots    \\ 
0 & B(t_s)- B(t_1) \\
\end{array}
\right)_{t=z}
\\
&=&
\begin{array}{c@{\,\,\,\,\,}c@{\,\,\,\,}|@{\,\,\,\,}c@{\,\,\,\,}@{\,\,\,\,}c@{\,\,\,\,}c
@{\,\,\,\,}c@{\,\,\,\,}c@{\,\,\,\,}ccc}
\ldelim({17}{4pt}[] && & &  &&\rdelim){17}{4pt}[]\\
&E_\ell &&\bigzerol & & \\
&&     &    &    \\ 
 \cline{2-6}  &&&&&&&&\rdelim\}{6}{10pt}[$\ell$ rows]\\
\\[-7mm]
&&\overrightarrow{\widetilde{f}(t_1)\widetilde{f}(t_2)} & & &  \bigzerol  & && \\   
&\bigzerol&&\overrightarrow{\widetilde{f}(t_1)\widetilde{f}(t_2)}  & \\
&&\bigzerol& &\ddots   &    \\ 
&&  & &&  \overrightarrow{\widetilde{f}(t_1)\widetilde{f}(t_2)}   \\
&&  & && \\
\\[-8mm]
 \cline{2-6} 
&\vdots& \vdots &\vdots &\vdots& \vdots  \\
  \cline{2-6}  &&&&&&&&\rdelim\}{6}{10pt}[$\ell$ rows]\\
\\[-7mm]
 &&  \overrightarrow{\widetilde{f}(t_1)\widetilde{f}(t_s)}   &  & &   \bigzerol       \\
&\bigzerol&& \overrightarrow{\widetilde{f}(t_1)\widetilde{f}(t_s)}   & &\\
&&\bigzerol&   &  \ddots   &    \\ 
&&  & &&   \overrightarrow{\widetilde{f}(t_1)\widetilde{f}(t_s)}    \\
\end{array}, 
\end{eqnarray*}
where 
$\overrightarrow{\widetilde{f}(t_1)\widetilde{f}(t_i)} =(
\widetilde{f}_1(t_i)-\widetilde{f}_1(t_1), 
\widetilde{f}_2(t_i)-\widetilde{f}_2(t_1), 
\ldots ,\widetilde{f}_m(t_i)-\widetilde{f}_m(t_1) 
)$ $(2\leq i \leq s)$ and $t=z$. 
From  $s-1\leq s_f-1$ and the definition of $s_f$, we have  
\begin{eqnarray*}
{\rm dim} \sum_{i=2}^s\mathbb{R}\overrightarrow{\widetilde{f}(t_1)\widetilde{f}(t_i)}=s-1,
\end{eqnarray*}
where $t=z$. Hence, 
by the construction of 
$Q_1M_2Q_2$ and $s-1\leq m$, the rank of 
$Q_1M_2Q_2$ 
is equal to $\ell s$. 
Therefore, the rank of $M_2$ 
must be equal to $\ell s$. 
Hence, we get (4.1). 
Therefore, there exists a subset $\Sigma$ of $\mathcal{L}(\mathbb{R}^m, \mathbb{R}^\ell)$ with Lebesgue measure zero such that 
for any $\pi \in \mathcal{L}(\mathbb{R}^m, \mathbb{R}^\ell)-\Sigma $, 
and  
for any $s$ $(2\leq s \leq s_f)$, 
the $C^r$ mapping 
$(F_\pi \circ f)^{(s)}:N^{(s)}\to (\R^\ell)^s$ is transverse to $\Delta_s$. 

Moreover, suppose that 
the $C^r$ mapping $F_\pi$ 
satisfies that $| F_\pi^{-1}(y) | \leq s_f$ for any $y\in\mathbb{R}^\ell$. 
Since $f:N\to \mathbb{R}^m$ is injective, 
it follows that $| (F_\pi \circ f)^{-1}(y) | \leq s_f$ for any $y\in\mathbb{R}^\ell$. 
Thus, for any positive integer $s$ with $s\geq s_f+1$, we have 
$(F_\pi \circ f)^{(s)}(N^{(s)})\bigcap \Delta_s=\emptyset$.  
Namely, for any positive integer $s$ with $s\geq s_f+1$, 
the $C^r$ mapping $(F_\pi \circ f)^{(s)}$ is transverse to $\Delta_s$. 
Hence, $F_\pi \circ f:N\to \mathbb{R}^\ell$ is a $C^r$ mapping with normal crossings. 
\hfill\qed
\section{Applications of \cref{main,main2}}\label{section 5}
In \cref{application1} (resp., \cref{application2}), 
applications of \cref{main} (resp., \cref{main2}) 
are stated and proved. 
In \cref{application2}, applications obtained by combining \cref{main,main2}   
are also given. 
\subsection{Applications of \cref{main}}\label{application1}

A $C^2$ function $g:N\to \mathbb{R}$ is called a {\it Morse function} 
if all of the critical points of $g$ are nondegenerate, 
where $N$ is a $C^2$ manifold of dimension $n$ 
(for details on Morse functions, see for instance, \cite[p. 63]{GG}). 
In the case of $(n,\ell)=(n,1)$, we have the following. 

\begin{corollary}\label{Morse function}
Let $f$ be a $C^2$ immersion 
of $N$ into an open subset $U$ of $\mathbb{R}^m$, 
where $N$ is a $C^2$ manifold of dimension $n$. 
Let $F:U\to \mathbb{R}$ be a $C^2$ function. 
Then, there exists a subset $\Sigma$ of $\mathcal{L}(\mathbb{R}^{m},\mathbb{R})$  with Lebesgue measure zero 
such that for any $\pi \in \mathcal{L}(\mathbb{R}^{m},\mathbb{R})-\Sigma$, 
the $C^2$ function $F_\pi \circ f:
N\to \mathbb{R}$ is a Morse function. 
\end{corollary}
{\it Proof.}\qquad 
We have $\dim N-\codim \Sigma^1(N,\R)=0$. 
Therefore, from \cref{main},  
there exists a subset $\Sigma$ with Lebesgue measure zero 
of $\mathcal{L}(\mathbb{R}^{m},\mathbb{R})$ 
such that for any $\pi \in \mathcal{L}(\mathbb{R}^{m},\mathbb{R})-\Sigma$, 
the mapping $j^1(F_\pi \circ f):
N\to J^1(N, \mathbb{R})$ is transverse to $\Sigma ^1(N,\mathbb{R})$. 
Therefore, if $q\in N$ is a critical point of the function $F_\pi \circ f$, then 
the point $q$ is nondegenerate. 
\hfill\qed
\par 
\bigskip

In the case of $\ell \geq2n$, we have the following. 
\begin{corollary}\label{immersion}
Let $f$ be a $C^2$ immersion 
of $N$ into an open subset $U$ of $\mathbb{R}^{m}$, 
where $N$ is a $C^2$ manifold of dimension $n$. 
Let $F:U\to \mathbb{R}^{\ell}$ be a $C^2$ mapping $(\ell \geq2n)$. 
Then, there exists a subset $\Sigma$ of 
$\mathcal{L}(\mathbb{R}^{m},\mathbb{R}^{\ell})$ with Lebesgue measure zero 
such that for any $\pi \in \mathcal{L}(\mathbb{R}^{m},\mathbb{R}^{\ell})-\Sigma$, 
the mapping $F_\pi \circ f:
N\to \mathbb{R}^{\ell}$ is a $C^2$ immersion. 
\end{corollary}
{\it Proof.}\qquad 
It is clearly seen that $F_\pi \circ f:N\to \mathbb{R}^\ell$ 
is an immersion 
if and only if   
$j^1(F_\pi \circ f)(N)\bigcap \bigcup_{k=1}^n\Sigma^k(N,\mathbb{R}^\ell)=
 \emptyset$. 
From $\ell \geq 2n$, 
for any positive integer $k$ $(1\leq k \leq n)$, we have 
\begin{eqnarray*}
\dim N-\codim \Sigma^k(N,\R^\ell) =n-k(\ell-n+k)\leq 0. 
\end{eqnarray*}
Thus, for any positive integer $k$ $(1\leq k \leq n)$, 
from \cref{main}, 
there exists a subset $\wt{\Sigma}_k$ of $\mathcal{L}(\mathbb{R}^{m},\mathbb{R}^{\ell})$ with Lebesgue measure zero 
such that for any 
$\pi \in \mathcal{L}(\mathbb{R}^{m},\mathbb{R}^{\ell})-\wt{\Sigma}_k$, 
the mapping $j^1(F_\pi \circ f):
N\to J^1(N,\mathbb{R}^{\ell})$ 
is transverse to $\Sigma^k(N,\mathbb{R}^\ell)$. 
Set $\Sigma=\bigcup_{k=1}^n\wt{\Sigma}_k$. 
Note that $\Sigma$ has Lebesgue measure zero.  
Let $\pi \in \mathcal{L}(\mathbb{R}^{m},\mathbb{R}^{\ell})-\Sigma$ be 
an arbitrary element.  
Then, suppose that there exists a point $q\in N$ 
and a positive integer $k$ $(1\leq k \leq n)$ 
such that $j^1(F_\pi \circ f)(q)\in \Sigma^k(N,\mathbb{R}^\ell)$. 
Since $j^1(F_\pi\circ f)$ is transverse to $\Sigma^k(N,\R^\ell)$, 
we have the following: 
\begin{eqnarray*}
d(j^1(F_\pi\circ f))_q(T_qN)+T_{j^1(F_\pi\circ f)(q)} \Sigma^k(N,\R^\ell)
=T_{j^1(F_\pi\circ f)(q)}J^1(N,\R^\ell). 
\end{eqnarray*}
Hence, we have 
{\small 
\begin{eqnarray*}
\dim d(j^1(F_\pi \circ f))_q(T_qN)
&\geq &
\dim T_{j^1(F_\pi \circ f)(q)}J^1(N,\R^\ell)-
\dim T_{j^1(F_\pi \circ f)(q)}\Sigma^k(N,\R^\ell)
\\
&=&
\codim  T_{j^1(F_\pi \circ f)(q)}\Sigma^k(N,\R^\ell). 
\end{eqnarray*}
}Thus, we get $n\geq k(\ell-n+k)$. 
This contradicts the assumption $\ell \geq 2n$. 
Therefore, we get 
$j^1(F_\pi\circ f)(N)\bigcap \bigcup_{k=1}^n\Sigma^k(N,\R^\ell)=\emptyset$. 
\hfill\qed
\par 
\bigskip 
A $C^1$ mapping $g:N\to \mathbb{R}^\ell$ 
{\it has singular points of corank at most k} 
if  
\begin{eqnarray*}
\sup\left\{{\rm corank\ }dg_q\mid q\in N \right\}\leq k, 
\end{eqnarray*}
where ${\rm corank\ }dg_q={\rm min} \{n,\ell \}-{\rm rank\ } dg_q$. 
\begin{corollary}\label{corank}
Let $f$ be a $C^r$ immersion 
of $N$ into an open subset $U$ of $\mathbb{R}^m$, 
where $N$ is a $C^r$ manifold of dimension $n$. 
Let $F:U\to \mathbb{R}^\ell$ be a $C^r$ mapping. 
Let $k_0$ be the maximum integer 
satisfying 
$(n-v+k_0)(\ell-v+k_0)\leq n$ $(v={\rm min}\{n,\ell\})$. 
If 
\[
r>\max \{ \dim N -\codim \Sigma ^1(N,\mathbb{R}^\ell),0\}+1, 
\]
then there exists a subset $\Sigma$ of $\mathcal{L}(\mathbb{R}^{m},\mathbb{R}^{\ell})$ 
with Lebesgue measure zero  
such that for any $\pi \in \mathcal{L}(\mathbb{R}^{m},\mathbb{R}^{\ell})-\Sigma$, 
the $C^r$ mapping $F_\pi \circ f:
N\to \mathbb{R}^\ell$ has 
singular points of corank at most $k_0$. 
\end{corollary}
{\it Proof.}\qquad 
For any positive integer $k$ $(1\leq k \leq v)$, 
we have 
\begin{eqnarray*}
r&>&\max \{ \dim N -\codim \Sigma ^1(N,\mathbb{R}^\ell),0\}+1
\\
&\geq &\max \{ \dim N -\codim \Sigma ^k(N,\mathbb{R}^\ell),0\}+1. 
\end{eqnarray*}

From \cref{main}, for any positive integer $k$ satisfying $1\leq k\leq v$, 
there exists a subset $\widetilde{\Sigma}_k$ of $\mathcal{L}(\mathbb{R}^{m},\mathbb{R}^{\ell})$ 
with Lebesgue measure zero 
such that for any $\pi \in \mathcal{L}(\mathbb{R}^{m},\mathbb{R}^{\ell})-\widetilde{\Sigma}_k$, 
the mapping $j^1(F_\pi \circ f):
N\to J^1(N,\mathbb{R}^{\ell})$ is transverse to $\Sigma^k(N,\mathbb{R}^\ell)$. 
Then, $\Sigma=\bigcup_{k=1}^v\widetilde{\Sigma}_k$ has Lebesgue measure zero. 
Hence, there exists a subset $\Sigma$ of $\mathcal{L}(\mathbb{R}^{m},\mathbb{R}^{\ell})$ 
with Lebesgue measure zero  
such that for any $\pi \in \mathcal{L}(\mathbb{R}^{m},\mathbb{R}^{\ell})-\Sigma$, 
the mapping $j^1(F_\pi \circ f):
N\to J^1(N,\mathbb{R}^{\ell})$ is transverse to 
$\Sigma ^k(N,\mathbb{R}^\ell)$ for any positive integer $k$ 
satisfying $1\leq k\leq v$. 

In the case of $\ell=1$, we have $k_0=1$. Thus, in this case, the assertion clearly holds. 

Now, we will consider the case of $\ell \geq 2$. 
In this case, note that $k_0+1\leq v$. Indeed, 
suppose that $v\leq k_0$. 
Then, by $(n-v+k_0)(\ell-v+k_0)\leq n$, we get $n\ell \leq n$. 
This contradicts the assumption $\ell \geq 2$. 
For the proof of \cref{corank}, it is sufficient to show that 
the mapping $j^1(F_\pi \circ f): 
N\to J^1(N,\mathbb{R}^{\ell})$ satisfies that 
$j^1(F_\pi \circ f)(N)\bigcap \Sigma^k(N,\mathbb{R}^\ell)=\emptyset$ 
for any positive integer $k$ satisfying $k_0+1\leq k\leq v$. 
Suppose that there exist a positive integer $k$ $(k_0+1\leq k \leq v)$ 
and a point $q\in N$ such that 
$j^1(F_\pi \circ f)(q)\in \Sigma^k(N,\mathbb{R}^\ell)$. 
Since the mapping $j^1(F_\pi \circ f):
N\to J^1(N,\mathbb{R}^{\ell})$ is transverse to $\Sigma^k(N,\mathbb{R}^\ell)$ 
at the point $q$, the following holds: 
\begin{eqnarray*}
d(j^1(F_{\pi}\circ f))_{q}(T_{q}N)
+T_{j^1(F_{\pi}\circ f)(q)}\Sigma ^{k}(N,\mathbb{R}^{\ell})
=T_{j^1(F_{\pi}\circ f)(q)}J^1(N,\mathbb{R}^{\ell}).
\end{eqnarray*}
Hence, we have
\begin{eqnarray*}
{}&&{}{\rm dim}\ d(j^1(F_{\pi}\circ f))_{q}(T_{q}N)\\
&\geq &{\rm dim}\ T_{j^1(F_{\pi}\circ f)(q)}J^1(N,\mathbb{R}^{\ell})-
{\rm dim}\ T_{j^1(F_{\pi}\circ f)(q)}\Sigma ^{k}(N,\mathbb{R}^{\ell})\\
&=&{\rm codim}\ T_{j^1(F_{\pi}\circ f)(q)}\Sigma ^{k}(N,\mathbb{R}^{\ell}).
\end{eqnarray*}
Thus, we get $n\geq (n-v+k)(\ell-v+k)$. 
Since the given integer $k_0$ is the maximum integer 
satisfying 
$n\geq(n-v+k_0)(\ell-v+k_0)$, 
it follows that $k\leq k_0$. 
This contradicts the assumption $k_0+1\leq k$. 
\hfill\qed
\par 
\bigskip 

\subsection{Applications of \cref{main2}}\label{application2}
\begin{corollary}\label{normal}
Let $f$ be a $C^r$ injection of $N$ into an open subset  
$U$ of $\mathbb{R}^m$, 
where $N$ is a $C^r$ manifold of dimension $n$. 
Let $F:U\to \mathbb{R}^\ell$ be a $C^r$ mapping. 
If 
\[
(s_f-1)\ell>n s_f \textit{ and } r>\max \{2n-\ell ,0\}, 
\]
then 
there exists a subset $\Sigma$ of $\mathcal{L}(\mathbb{R}^m, \mathbb{R}^\ell)$ 
with Lebesgue measure zero such that 
for any $\pi \in \mathcal{L}(\mathbb{R}^m, \mathbb{R}^\ell)-\Sigma $, 
$F_\pi \circ f:N\to \mathbb{R}^\ell$ is a $C^r$ mapping with normal crossings 
satisfying $(F_\pi \circ f)^{(s_f)}(N^{(s_f)})\bigcap\Delta _{s_f}=\emptyset$. 
\end{corollary}
{\it Proof.}\qquad 
From $(s_f-1)\ell>n s_f$, we have $n-\ell<0$. Thus, we get 
\begin{eqnarray*}
s_0&=&\max \{ s(n-\ell)+\ell \ |\ 2 \leq s \leq s_f\}
\\
&=&2n-\ell. 
\end{eqnarray*}
Hence, note that $r>\max\{s_0,0\}$. 
From \cref{main2}, 
there exists a subset $\Sigma$ of $\mathcal{L}(\mathbb{R}^m, \mathbb{R}^\ell)$ with Lebesgue measure zero such that 
for any $\pi \in \mathcal{L}(\mathbb{R}^m, \mathbb{R}^\ell)-\Sigma $, 
and  
for any $s$ $(2\leq s \leq s_f)$, 
the mapping $(F_\pi \circ f)^{(s)}:N^{(s)}\to (\R^\ell)^s$ is transverse to $\Delta_s$. 
Therefore, for this proof,  
it is sufficient to prove that 
for any $\pi \in \mathcal{L}(\mathbb{R}^m, \mathbb{R}^\ell)-\Sigma $, 
the mapping $(F_\pi \circ f)^{(s_f)}$ satisfies that 
$(F_\pi \circ f)^{(s_f)}(N^{(s_f)})\bigcap \Delta_{s_f}=\emptyset$. 

Suppose that there exists an element 
$\pi \in \mathcal{L}(\mathbb{R}^m, \mathbb{R}^\ell)-\Sigma $ 
such that there exists a point $q\in N^{(s_f)}$ 
satisfying $(F_{\pi}\circ f)^{(s_f)}(q)\in \Delta_{s_f}$. 
Since $(F_{\pi}\circ f)^{(s_f)}$ is transverse to $\Delta_{s_f}$, 
we have the following: 
\begin{eqnarray*}
d((F_{\pi}\circ f)^{(s_f)})_{q}(T_{q}N^{(s_f)})
+T_{(F_{\pi}\circ f)^{(s_f)}(q)}\Delta_{s_f}
=T_{(F_{\pi}\circ f)^{(s_f)}(q)}(\mathbb{R}^\ell)^{s_f}.
\end{eqnarray*}
Thus, we get 
\begin{eqnarray*}
{}&&{}{\rm dim}\ d((F_{\pi}\circ f)^{(s_f)})_{q}(T_{q}N^{(s_f)})\\
&\geq &{\rm dim}\ T_{(F_{\pi}\circ f)^{(s_f)}(q)}(\mathbb{R}^\ell)^{s_f}-
{\rm dim}\ T_{(F_{\pi}\circ f)^{(s_f)}(q)}\Delta_{s_f}\\
&=&{\rm codim}\ T_{(F_{\pi}\circ f)^{(s_{f})}(q)}\Delta_{s_f}.
\end{eqnarray*}
Hence, we have $ns_f\geq (s_f-1)\ell$. 
This contradicts the assumption $(s_f-1)\ell>ns_f$. 
\hfill\qed
\par 
\bigskip  
In the case of $\ell >2n$, 
we have the following. 
\begin{corollary}\label{injective}
Let $f$ be a $C^1$ injection of $N$ into an open subset  
$U$ of $\mathbb{R}^m$, 
where $N$ is a $C^1$ manifold of dimension $n$. 
Let $F:U\to \mathbb{R}^\ell$ be a $C^1$ mapping. 
If $\ell>2n$, then 
there exists a subset $\Sigma$ of $\mathcal{L}(\mathbb{R}^m, \mathbb{R}^\ell)$ 
with Lebesgue measure zero such that 
for any $\pi \in \mathcal{L}(\mathbb{R}^m, \mathbb{R}^\ell)-\Sigma $, 
the $C^1$ mapping $F_\pi \circ f:N\to \mathbb{R}^\ell$ is injective. 
\end{corollary}
{\it Proof.}\qquad 
Since $s_f\geq 2$ and $\ell>2n$, it is easily seen that the dimension pair $(n, \ell)$ 
satisfies the assumption $(s_f-1)\ell>ns_f$ of \cref{normal}. 
Indeed, from $\ell>2n$, we get $(s_f-1)\ell>2n(s_f-1)$. 
From $s_f\geq 2$, it follows that $2n(s_f-1)\geq ns_f$. 
 
Since $\max\{2n-\ell,0\}=0$, 
from \cref{normal}, 
there exists a subset $\Sigma$ of $\mathcal{L}(\mathbb{R}^m, \mathbb{R}^\ell)$ 
with Lebesgue measure zero such that 
for any $\pi \in \mathcal{L}(\mathbb{R}^m, \mathbb{R}^\ell)-\Sigma $, 
the mapping $(F_\pi \circ f)^{(2)} : N^{(2)} \to (\mathbb{R}^\ell)^2$ is 
transverse to $\Delta_2$. 
For this proof, 
it is sufficient to prove that 
the mapping $(F_\pi \circ f)^{(2)}$ satisfies that 
$(F_\pi \circ f)^{(2)}(N^{(2)})\bigcap \Delta_2 =\emptyset$. 

Suppose that there exists a point $q\in N^{(2)}$ such that 
$(F_\pi \circ f)^{(2)}(q)\in \Delta_2$. 
Then, we get the following: 
\begin{eqnarray*}
d((F_{\pi}\circ f)^{(2)})_{q}(T_{q}N^{(2)})
+T_{(F_{\pi}\circ f)^{(2)}(q)}\Delta_{2}
=T_{(F_{\pi}\circ f)^{(2)}(q)}(\mathbb{R}^\ell)^{2}.
\end{eqnarray*}
Thus, we have
\begin{eqnarray*}
{}&&{}{\rm dim}\ d((F_{\pi}\circ f)^{(2)})_{q}(T_{q}N^{(2)})\\
&\geq &{\rm dim}\ T_{(F_{\pi}\circ f)^{(2)}(q)}(\mathbb{R}^\ell)^{2}-
{\rm dim}\ T_{(F_{\pi}\circ f)^{(2)}(q)}\Delta_{2}\\
&=&{\rm codim}\ T_{(F_{\pi}\circ f)^{(2)}(q)}\Delta_{2}.
\end{eqnarray*}
Hence, we have $2n\geq \ell$. 
This contradicts the assumption $\ell>2n$. 
\hfill\qed
\par 
\bigskip  

By combining \cref{immersion,injective}, 
we have the following. 
\begin{corollary}\label{injective immersion}
Let $f$ be an injective immersion of $N$ into an open subset  
$U$ of $\mathbb{R}^m$, 
where $N$ is a $C^2$ manifold of dimension $n$ and 
$f$ is of class $C^2$. 
Let $F:U\to \mathbb{R}^\ell$ be a $C^2$ mapping. 
If $\ell>2n$, then 
there exists a subset $\Sigma$ of $\mathcal{L}(\mathbb{R}^m, \mathbb{R}^\ell)$ 
with Lebesgue measure zero such that 
for any $\pi \in \mathcal{L}(\mathbb{R}^m, \mathbb{R}^\ell)-\Sigma $, 
the $C^2$ mapping $F_\pi \circ f:N\to \mathbb{R}^\ell$ is an injective immersion. 
\end{corollary}

From \cref{injective immersion}, we get the following. 
\begin{corollary}\label{embedding}
Let $N$ be a compact $C^2$ manifold of dimension $n$. 
Let $f$ be a $C^2$ embedding of $N$ into an open subset  
$U$ of $\mathbb{R}^m$. 
Let $F:U\to \mathbb{R}^\ell$ be a $C^2$ mapping. 
If $\ell>2n$, then 
there exists a subset $\Sigma$ of $\mathcal{L}(\mathbb{R}^m, \mathbb{R}^\ell)$ 
with Lebesgue measure zero such that 
for any $\pi \in \mathcal{L}(\mathbb{R}^m, \mathbb{R}^\ell)-\Sigma $, 
the $C^2$ mapping $F_\pi \circ f:N\to \mathbb{R}^\ell$ is an embedding. 
\end{corollary}
\section{Proof of \cref{abra}}\label{app}
\subsection{Preliminaries for the proof of \cref{abra}}
Let $N$ and $P$ be $C^r$ manifolds, and let $g : N\to P$ be a $C^1$ mapping 
($r\geq1$). A point $x\in N$ is called a {\it critical point} of $g$ if
it is not a regular point, i.e., the rank of $dg_x$ is less than the dimension of
$P$. We say that a point $y\in P$ is a {\it critical value} if it is the image of a critical point. 
A point $y\in P$ is called a {\it  regular value} if it is not a critical value. 
The following is Sard's theorem. 
\begin{theorem}[\cite{sard}]\label{sard}
If $N$ and $P$ are $C^r$ manifolds, 
$g : N \to P$ is a $C^r$ mapping, and $r> \max\{\dim N-\dim P,0\}$, then the set of
critical values of $g$ has Lebesgue measure zero.
\end{theorem}

\subsection{Proof of \cref{abra}}
In this proof, by $\pi:N\times A\to A$, we denote 
the natural projection 
defined by $\pi(x,a)=a$. 

Since $\Gamma$ is transverse to $Z$, the set $\Gamma^{-1}(Z)$ is 
a $C^r$ submanifold of $N\times A$ satisfying 
\begin{eqnarray}\label{eq:lemma1r}
\dim N+\dim A- \dim \Gamma^{-1}(Z)=\dim P-\dim Z. 
\end{eqnarray}

Firstly, suppose that $\dim \Gamma^{-1}(Z)=0$. Then, 
since $\Gamma^{-1}(Z)$ is a countable set, 
$\pi(\Gamma^{-1}(Z))$ has Lebesgue measure zero in $A$. 
It is clearly seen that for any $a\in A-\pi(\Gamma^{-1}(Z))$, 
the mapping $\Gamma_a$ is transverse to $Z$. 

Finally, we will consider the case $\dim \Gamma^{-1}(Z)>0$. 
It is not hard to see that if $a\in A$ is 
a regular value of $\pi|_{\Gamma^{-1}(Z)}$, then 
$\Gamma_a$ is transverse to $Z$, 
where  $\pi|_{\Gamma^{-1}(Z)}$ is the restriction of $\pi$ to $\Gamma^{-1}(Z)$. 
Let $\Sigma$ be the set of critical values of $\pi|_{\Gamma^{-1}(Z)}$. 
From $r>\max \{ \dim N+\dim Z-\dim P,0\}$ and 
(\ref{eq:lemma1r}), we have $r>\max \{ \dim \Gamma^{-1}(Z)-\dim A,0\}$. 
From \cref{sard}, $\Sigma$ has Lebesgue measure zero in $A$. 
Therefore, if $a\in A-\Sigma$, then 
$\Gamma_a$ is transverse to $Z$. 
\hfill\qed

\section*{Acknowledgements}
The author is most grateful to the anonymous reviewer for carefully reading 
the first manuscript of this paper and for giving invaluable suggestions. 
The author was supported by JSPS KAKENHI Grant Number 16J06911. 
\>

\end{document}